\documentclass[a4paper, 11pt]{article}

\usepackage{amsmath,amsfonts,amssymb,amsthm}
\usepackage[a4paper, margin=1in]{geometry}
\usepackage{enumitem}

\usepackage[T1]{fontenc}

\usepackage{lmodern}
\usepackage{calc}
\usepackage[dvipsnames]{xcolor}
\usepackage{nameref}
\usepackage{upgreek}     
\usepackage{enumitem}
\usepackage{graphicx}
\usepackage{mathrsfs}
\usepackage{tikz-cd}
\usepackage{hhline}
\usepackage{textcomp}
\usepackage{amscd}
\usepackage[all,cmtip]{xy}
\usepackage{mathtools}
\usepackage{extarrows}
\usepackage{upgreek}
\usepackage{bbm}
\usepackage[displaymath]{lineno}
\usepackage[singlespacing]{setspace}%
\usepackage[colorlinks=true,breaklinks=true,bookmarks=true,urlcolor=blue,
citecolor=BrickRed,linkcolor=teal,bookmarksopen=false,draft=false]{hyperref}

\providecommand{\U}[1]{\protect\rule{.1in}{.1in}}

\makeatletter
\let\orgdescriptionlabel\descriptionlabel
\renewcommand*{\descriptionlabel}[1]{%
	\let\orglabel\label
	\let\label\@gobble
	\phantomsection
	\edef\@currentlabel{#1}%
	\let\label\orglabel
	\orgdescriptionlabel{#1}%
}
\makeatother

\theoremstyle{plain}
\newtheorem{theorem}{Theorem}[section]
\newtheorem{lemma}[theorem]{Lemma}
\newtheorem{corollary}[theorem]{Corollary}

\theoremstyle{definition}
\newtheorem{definition}[theorem]{Definition}
\newtheorem*{definition*}{Definition} 
\newcommand{\bl}[1]{\textup{#1}}
\newcommand{\RR}{\mathbb{R}}

\newcommand{\ps}{\mathrm{({PS}})_m}

\sbox0{$\fontdimen8\textfont3=0.4pt$}

\newcommand{\abs}[1]{\lvert#1\rvert}

\DeclareMathAlphabet{\mathpzc}{OT1}{pzc}{m}{it}

\newcommand{\ftt}[1] {\mathsf{#1}}

\newcommand{\va}{\varphi}

\newcommand{\dd}{{\tt D}}

\newcommand{\fs}[1]{\mathsf {#1}}

\newcommand{\mt}{\mathbbm {d}}

\newcommand{\set}[1]{\left\{#1\right\}}
\newcommand\Set[2]{\left\{#1\mid#2\right\}} 
\newcommand{\snorm}[2][]{\left\lVert#2\right\rVert_{#1}}

\newcommand{\Sem}[1]  {\textsf{Sem}(#1)}

\newcommand{\zero}[1]{\boldsymbol{0}_{#1}}

\usepackage{bm}    

\newcommand{\rr}{\mathbb{R}}
\newcommand{\nn}{\mathbb{N}}

\newcommand{\bb}{\mathfrak{S}}
\DeclareMathAlphabet\EuScript{U}{eus}{m}{n}
\SetMathAlphabet\EuScript{bold}{U}{eus}{b}{n}
\newcommand\opn{\ensuremath{\mathrel{\mathpalette\opncls\circ}}}

\newcommand{\opncls}[2]{
	\ooalign{$#1\subseteq$\cr
		\hidewidth\raisefix{#1}\hbox{$#1{\stylefix{#1}#2}\mkern2mu$}\cr}}
\def\raisefix#1{
	\ifx#1\displaystyle
	\raise.39ex
	\else
	\ifx#1\textstyle
	\raise.39ex
	\else
	\ifx#1\scriptstyle
	\raise.275ex
	\else
	\raise.150ex
	\fi
	\fi
	\fi
}
\def\stylefix#1{
	\ifx#1\displaystyle
	\scriptstyle
	\else
	\ifx#1\textstyle
	\scriptstyle
	\else
	\ifx#1\scriptstyle
	\scriptscriptstyle
	\else
	\scriptscriptstyle
	\fi
	\fi
	\fi
}
\DeclareFontFamily{U}{mathx}{\hyphenchar\font45}
\DeclareFontShape{U}{mathx}{m}{n}{
	<5> <6> <7> <8> <9> <10>
	<10.95> <12> <14.4> <17.28> <20.74> <24.88>
	mathx10
}{}
\newcommand{\fr}{Fr\'{e}chet }

\DeclareMathAlphabet{\mathsfit}{OT1}{cmss}{m}{sl}

\newcommand{\subjclass}[1]{\textbf{AMS Subject Classifications (2010):} #1\par}
\newcommand{\keywords}[1]{\textbf{Keywords:} #1\par}
\title{A Generalized Palais-Smale Condition in the Fr\'{e}chet space setting}
\author{Kaveh Eftekharinasab}
\date{In memory of Vladimir Sharko}
\usepackage{cite}
\begin{document}
	
	\maketitle
	
	\begin{abstract}
We extend the Palais-Smale condition to Keller's $C_c^1$-functionals on Fr\'{e}chet spaces. Using this condition together with  Ekeland's variational principle, we obtain some results regarding the existence of minima. In this setting, we prove that the Palais-Smale condition for functionals bounded below
implies the coercivity.
	\end{abstract}
\bigskip

\paragraph{Erratum.}	
In the published version of this paper, we defined the Palais-Smale condition on \fr manifolds using the concept of the $\beta$-cotangent bundle, as introduced in \cite{tw}. However, this definition has since been shown to be invalid (see \cite[Remark 9.7]{gl}). Therefore, in this updated version, we have removed the Palais-Smale definition and all related results concerning \fr manifolds, and instead refer the reader to \cite{eftekharinasab}.
	\let\thefootnote\relax\footnotetext{
		\,\,\,\,\,\keywords{The Palais-Smale condition, Fr\'{e}chet spaces, Keller's $C_c^1$-functionals, Ekelend's variational principle, coercivity.}}
	\section*{Introduction}
The Palais-Smale  condition is a fundamental tool in critical point theory for infinite-dimensional spaces. Since its introduction by Palais and Smale~\cite{ps}, numerous generalizations have been developed for various classes of maps (see~\cite{mm} for a survey). These conditions have been crucial in addressing numerous problems within critical point theory and its applications across differential equations, geometry, and physics (see~\cite{mm}).

In this paper, we extend the Palais-Smale condition to Keller's $C_c^1$-functionals on Fr\'{e}chet spaces. It is well-established that deformation lemmas are essential tools for locating critical points of $C^1$-functionals. These deformations typically arise as solutions to Cauchy problems. However, on general Fr\'{e}chet spaces, ordinary differential equations may lack unique solutions for a given initial condition, hindering the direct application of deformation lemma techniques. Nevertheless, we can still employ Ekeland's variational principle to locate critical points. We apply the version of Ekeland's variational principle for Fr\'{e}chet spaces due to Qui~\cite{ekel} to obtain standard results concerning the existence of minima.

Furthermore, we add a new result to the published version: that the Palais-Smale condition for functionals bounded below implies coercivity.
\section{The Palais-Smale condition}

Throughout this paper, we assume that \( (\fs{F}, \Sem{\fs{F}}) \) and \( (\fs{E}, \Sem{\fs{E}}) \) are \fr spaces over \( \rr \), where \( \Sem{\fs{F}} = \Set{\snorm[\fs{F},n]{\cdot}}{n \in \nn} \) and \( \Sem{\fs{E}} = \Set{\snorm[\fs{E},n]{\cdot}}{n \in \nn} \) are increasing sequences of continuous seminorms that define the topologies of \( \fs{F} \) and \( \fs{E} \), respectively. 

We use the notation \( U \opn \mathsf{T} \) to denote that \( U \) is an open subset of the topological space \( \mathsf{T} \).

Let \( \mathfrak{S} \) be a family of bounded subsets of \( \fs{F} \), with the following two properties:
\begin{itemize}
	\item[(\(\mathfrak{S}_1\))] If \( A, B \in \mathfrak{S} \), then there exists \( C \in \mathfrak{S} \) such that \( A \cup B \subset C \).
	\item[(\(\mathfrak{S}_2\))] If \( A \in \mathfrak{S} \) and \( r \) is a real number, then there exists \( B \in \mathfrak{S} \) such that \( rA \subset B \).
\end{itemize}

Throughout the paper, we assume that \( \mathfrak{S} \) is the family of all compact subsets of \( \fs{F} \). We consider the topology of compact convergence on the dual space, and denote the dual space endowed with this topology with by \( \fs{F}'_c \).

 Let \( \mathcal{L}_c(\fs{F}, \fs{E}) \) be the space of all continuous linear mappings from \( \fs{F} \) to \( \fs{E} \) endowed with the topology of compact convergence, which is a Hausdorff locally convex topology, defined by seminorms:
 \[
 \snorm[{S},i]{\ell} \coloneqq \sup\Set{\snorm[\fs{E},i]{\ell (f)}}{ f  \in S}
\]
where \(S \in \mathfrak{S}\). If \(\fs{E}=\RR \) with the usual modulus \(\abs{\cdot}\), then
we denote by \(\snorm[{S}]{\cdot}\) the seminorms that define the topology of \( \mathcal{L}_c(\fs{F}, \rr) \) 
\begin{definition}[Definition 1.0.0, \cite{ke}]\label{def:diff}
	Let  $ \va\colon U \opn \fs{E}  \to  \fs{F}$ be a mapping. Then the  derivative
	of $\va$ at $x$ in the direction $h$ is defined by 
	\[
	\dd \va_x(h)=\dd\va(x)(h) \coloneqq
	\lim_{t \to 0} {1\over t}(\va(x+th) -\va(x))
	\]
	whenever it exists. 
	The function $\va$ is called differentiable at
	$x$ if $\dd \va(x)(h)$ exists for all $h \in \fs{E}$. 
	It is called \(C_c^1\)-mapping
	if it is differentiable at all
	points of $U$, and the mapping
	\[
	\dd \va \colon U  \to  \mathcal{L}_c(\fs{E}, \fs{F}) 
	\]
	is continuous. 
\end{definition}
 Higher order differentiability is defined in \cite[Definition 2.5.0]{ke}

The primary motivation for employing this class of mappings is the need to equip dual spaces with a suitable topology in order to define the Palais–Smale condition. These mappings are known as Keller \( C_c^k \)-mappings, which are equivalent to the well-established and widely used Michal--Bastiani notion of differentiability.

\begin{definition}[The Palais-Smale Condition]\label{def:PS}
	Let $\va\colon U \opn \fs{F} \to \rr$ be a Keller $C_c^1$-functional.
	\begin{itemize}
		\item [(i)] We say that $\va$ satisfies the Palais-Smale condition, $\mathrm{({PS}})$-condition in short, if every sequence $(x_i) \subset \fs{F}$ such that $\va(x_i)$ is bounded and
		$$
		\dd \va(x_i) \to 0 \quad \text{in} \quad \fs{F}_{c}',
		$$
		has a convergent subsequence.
		\item[(ii)] We say that $\va$ satisfies the Palais-Smale condition at level $m \in \rr$, $\mathrm{({PS}})_m$-condition in short, if every sequence $(x_i) \subset \fs{F}$ such that
		$$
		\va(x_i) \to m \quad \text{and} \quad \dd \va(x_i) \to 0 \quad \text{in} \quad \fs{F}_{c}',
		$$
		has a convergent subsequence.
	\end{itemize}
\end{definition}

If $ \va $ satisfies the $ \ps $-condition, then every $ \ps $-sequence converges, up to a subsequence,
to some point $ p $  and, by continuity, one has that  $ \va(p)=m $ and $ \va'(p) =0 $. In another word, $ p $
is a critical point of $ \va $. 

 We now show that the Palais-Smale condition provides a way to distinguish a certain family of neighborhoods of critical points and a means to characterize regular values of functionals.

Let \( \{ \Set{\snorm[S]{\cdot}}{S \in \mathfrak{S}}\} \) be a saturated family of seminorms defining the topology of the compact dual \( \fs{F}'_{c} \). Let \( m \in \rr \), \( r, s \in \rr_{+} \), \( n \in \nn \), and \( S \in \bb \). Define the following sets:

\begin{description}
	\item[{Critical set}] 
	\( \ftt{Cr}(\va) \coloneqq \Set{x \in \fs{F}}{\va'(x)=0} \);
	\item[{Critical level set}] 
	\( \ftt{Cr}(\va,m) \coloneqq \Set{x \in \ftt{Cr}(\va)}{\va(x) = m} \);
	\item[{\quad}] 
	\( \mathcal{U}_{S}(m,r) \coloneqq \Set{x \in \fs{F}}{\abs{\va(x)-m} < r,\, \snorm[S]{\va'(x)} < r} \);
	\item[{\quad}] 
	\( \mathcal{V}_{S,n}(m,s) \coloneqq \bigcup_{y \in \ftt{Cr}(\va,m)} \Set{x \in \fs{F}}{\snorm[\fs{F},n]{x - y} < s} \).
\end{description}

\begin{lemma}\label{lem:scs}
	Let \( \va: \fs{F} \to \rr \) be a $C_c^1$-functional satisfying the \( \ps \)-condition for some \( m \in \rr \). Then the following hold$\colon$
	\begin{enumerate}[label=\bl{(I.\arabic*)}, ref=I.\arabic*]
		\item \label{lem:psn1} The set \( \ftt{Cr}(\va,m) \) is compact.
		\item The family \( \Set{\mathcal{U}_{S}(m,r)}{S \in \bb,\, r \in \rr_{+}} \) is a fundamental system of neighborhoods of \( \ftt{Cr}(\va,m) \).
		\item The family \( \Set{\mathcal{V}_{S,n}(m,s)}{n \in \nn,\, S \in \bb, s \in \rr_{+}} \) is a fundamental system of neighborhoods of \( \ftt{Cr}(\va,m) \).
	\end{enumerate}
\end{lemma}

\begin{proof}
	\bl{(I.1)} By the \( \ps \)-condition, any sequence \( (x_i) \subset \ftt{Cr}(\va,m) \) has a convergent subsequence. The continuity of \( \va \) and \( \va' \) ensures that accumulation points lie in \( \ftt{Cr}(\va,m) \), hence it is compact.
	
	\medskip
	
	\bl{(I.2)} Each \( \mathcal{U}_{S}(m,r) \) with \( r \in \rr_{+}, S \in \bb_{\mathrm{k}} \) is a neighborhood of \( \ftt{Cr}(\va,m) \). Conversely, let \( \mathcal{N} \) be arbitrary open neighborhood of \( \ftt{Cr}(\va,m) \). Suppose for some sequence \( r_k \to 0 \), there exist points \( x_k \in \mathcal{U}_{S}(m,r_k) \setminus \mathcal{N} \). If \( (y_k) \subset \ftt{Cr}(\va,m) \) satisfies
	\[
	\snorm[\fs{F},n]{x_k - y_k} \leq r_k \quad \forall n \in \nn,
	\]
	then compactness of \( \ftt{Cr}(\va,m) \) implies \( y_k \to y \in \ftt{Cr}(\va,m) \), hence \( x_k \to y \in \mathcal{N} \), which contradicts the assumption.
	
	\medskip
	
	\bl{(I.3)} Each \( \mathcal{V}_{S,n}(m,s) \) with \( s \in \rr_{+}, S \in \bb_{\mathrm{k}}, n \in \nn \) is a neighborhood of \( \ftt{Cr}(\va,m) \). Suppose \( \mathcal{N} \) is an open neighborhood and for a sequence \( s_k \to 0 \), there exist points \( x_k \in \mathcal{V}_{S,n}(m,s_k) \setminus \mathcal{N} \). The \( \ps \)-condition implies \( (x_k) \) accumulates at some critical point \( u \in \ftt{Cr}(\va,m) \subset \mathcal{N} \), leading again to a contradiction.
\end{proof}

Although, Lemma \ref{lem:scs}\eqref{lem:psn1} states that the critical level sets of $ \va $, $ \ftt{Cr}(\va,m)$, are compact. However, satisfying the PS-condition at all levels does not imply that the critical set is bounded.

Now, we prove a minimization result for functionals that are bounded below by means of a version of the Ekeland variational principle for \fr spaces. 
\begin{theorem}[Corollary 2.1, \cite{ekel}] \label{th:ekef}
	Let a functional $f : \fs{F} \rightarrow (-\infty, \infty]$ be lower semi-continuous,  bounded
	from below and not  identically equal to $+\infty$.  Let $ \eta >0$ and $x_0 \in \fs{F}$ be given such that $f(x_0) \leqq \inf_{F}f + \eta$,
	and let $(\lambda_n)_{n \in \nn}$ be a sequence of positive real numbers. Then for any $i \in \nn$, there exists
	$z \in F$ such that
	\begin{enumerate}
		\item $\lambda_j \snorm[\fs{F},j]{z-x_0} \leq f(x_0) - f(z)$ for $j= 1,\ldots,i$;
		\item $  \snorm[\fs{F},j]{z-x_0} < \dfrac{\eta}{\lambda_j} $ for $j= 1,\ldots,i$;
		\item for any $x\in F, x \ne z $, there exists $m \in \nn$ such that 
		\[
		\lambda_m \snorm[\fs{F},m]{x-z} + f(x) > f(z),
		\]
		or equivalently,  $\sup_{k \in \nn} \lambda_n \snorm[\fs{F},m]{x-z} + f(x) > f(z)$ for any $x\in F, x \neq z $.
	\end{enumerate}
\end{theorem}
\begin{theorem} \label{sec}
	Let $\va\colon U \opn \fs{F}  \to \rr$ be a Keller $C_c^1$-functional.
Let $\va$ be bounded from below. 
Then, for each $\varepsilon>0$, $i \in \nn$, $x \in \fs{F}$ such that $\va(x) \leqq \inf_{\fs{F}}f + \varepsilon$,
there exist $z \in \fs{F}$ such that 
	\begin{enumerate}[label=\textup{(\alph*)},ref=\textup{\alph*}]
		\item \label{item:a} $\snorm[\fs{F},j]{z - x_0} \le \dfrac{\va(x_0) - \va(z)}{\sqrt{\varepsilon}}$ for $j = 1, \ldots, i$;
		\item \label{item:b} $\snorm[\fs{F},j]{z - x_0} < \sqrt{\varepsilon}$ for $j = 1, \ldots, i$;
		\item \label{item:c} $\snorm[S]{\dd \va(z)} \le K_S \sqrt{\varepsilon}$, for all $S \in \mathfrak{S}$ and some $K_S > 0$.
	\end{enumerate}
\end{theorem}

\begin{proof}
	Given $\varepsilon >0$ and $i \in \nn$. Let $x_0 \in \fs{F}$ be such that $\va(x_0) \le \inf_{\fs{F}}\va + \varepsilon$.
	Apply Theorem~\ref{th:ekef} with $\lambda_j = \sqrt{\varepsilon}$ for all $j \in \nn$ and $\eta = \varepsilon$.
	This yields a point $z \in \fs{F}$ (which we denote as $z$ instead of $z_{\varepsilon}$ for simplicity as $\varepsilon$ is fixed in the theorem statement) such that conditions (1) and (2) of Theorem \ref{th:ekef} hold:
	\begin{enumerate}
		\item $\sqrt{\varepsilon} \snorm[\fs{F},j]{z-x_0} \le \va(x_0) - \va(z)$ for $j= 1,\ldots,i$;
		\item $\snorm[\fs{F},j]{z-x_0} < \dfrac{\varepsilon}{\sqrt{\varepsilon}} = \sqrt{\varepsilon}$ for $j= 1,\ldots,i$.
	\end{enumerate}
	These correspond precisely to conditions \eqref{item:a} and \eqref{item:b} of the current theorem.
	For condition (3) of Theorem~\ref{th:ekef}, there exists $m_0 \in \nn$ such that for any $x \in \fs{F}, x \ne z$, we have
	\begin{equation}\label{eq:ekef_applied_specific}
		\sqrt{\varepsilon} \snorm[\fs{F},m_0]{x-z} + \va(x) > \va(z).
	\end{equation}
	Rearranging \eqref{eq:ekef_applied_specific}, we get
	\[
	\va(x) - \va(z) > -\sqrt{\varepsilon} \snorm[\fs{F},m_0]{x-z} \quad \forall x \ne z \in \fs{F}. 
	\]
	Now, let $h \in \fs{F}$ be arbitrary. For $t > 0$ sufficiently small such that $z+th \in U$ and $z+th \ne z$, we can replace $x$ with $z+th$:
	\[
	 \va(z+th) - \va(z) > -\sqrt{\varepsilon} \snorm[\fs{F},m_0]{th}.
	\]
	Dividing by $t$ (which is positive) gives
	\[
	\frac{\va(z+th) - \va(z)}{t} > -\sqrt{\varepsilon} \snorm[\fs{F},m_0]{h}. 
	\]
	Taking the limit as $t \to 0^+$ (since $\va$ is Keller $C_c^1$, it is Gateaux differentiable), we obtain
	$$ \dd \va(z)(h) \ge -\sqrt{\varepsilon} \snorm[\fs{F},m_0]{h} \quad \forall h \in \fs{F}. $$
	Now, replace $h$ with $-h$:
	$$ \dd \va(z)(-h) \ge -\sqrt{\varepsilon} \snorm[\fs{F},m_0]{-h}. $$
	Since $\dd \va(z)$ is linear and $\snorm[\fs{F},m_0]{\cdot}$ is a seminorm (hence $\snorm[\fs{F},m_0]{-h} = \snorm[\fs{F},m_0]{h}$), this becomes
	$$ -\dd \va(z)(h) \ge -\sqrt{\varepsilon} \snorm[\fs{F},m_0]{h}, $$
	which implies
	$$ \dd \va(z)(h) \le \sqrt{\varepsilon} \snorm[\fs{F},m_0]{h} \quad \forall h \in \fs{F}. $$
	Combining the two inequalities, we have
	$$ \abs{\dd \va(z)(h)} \le \sqrt{\varepsilon} \snorm[\fs{F},m_0]{h} \quad \forall h \in \fs{F}. $$
	Now consider $\snorm[S]{\dd \va(z)}$ for any $S \in \mathfrak{S}$:
	$$ \snorm[S]{\dd \va(z)} = \sup_{h \in S} \abs{\dd \va(z)(h)}. $$
	Since $S$ is a compact subset of $\fs{F}$, for each seminorm $\snorm[\fs{F},n]{\cdot}$, there exists a constant $C_{S,n}$ such that $\snorm[\fs{F},n]{h} \le C_{S,n}$ for all $h \in S$. Specifically, for the seminorm $\snorm[\fs{F},m_0]{\cdot}$ involved in the inequality:
	\[
	\snorm[S]{\dd \va(z)} = \sup_{h \in S} \abs{\dd \va(z)(h)} \le \sup_{h \in S} (\sqrt{\varepsilon} \snorm[\fs{F},m_0]{h}) = \sqrt{\varepsilon} \sup_{h \in S} \snorm[\fs{F},m_0]{h}. 
	\]
	Let $K_S = \sup_{h \in S} \snorm[\fs{F},m_0]{h}$. Since $S$ is compact (and thus bounded), $K_S$ is a finite positive constant for any $S \in \mathfrak{S}$ (if $S$ is non-empty).
	Therefore, we obtain
	 $\snorm[S]{\dd \va(z)} \le K_S \sqrt{\varepsilon} $.
\end{proof}

Let $\va : \fs{F} \to \rr$ be a functional bounded from below. A sequence $(x_i) \subset \fs{F}$ is said to be \textit{minimizing} if
$$ \lim_{i \to \infty}\va(x_i) = \inf_{x \in \fs{F}} \va(x). $$
A minimizing sequence $(x_i) \subset \fs{F}$ is said to consist of almost critical points of $\va$ if
$$ \dd \va(x_i) \to 0 \quad \text{in} \quad \fs{F}'_c. $$

\begin{corollary}\label{cor:almostc}
	Let $\va$ be as in Theorem~\ref{sec}. Then for any minimizing sequence $(x_i)_{i \in \nn}$, there exists a minimizing sequence $(y_i)_{i \in \nn}$ such that
	\begin{enumerate}[label=\textup{(\arabic*)},ref=\textup{\arabic*}]
		\item \label{eq:cms1} $\va(x_i) \le \va(y_i)$;
		\item \label{eq:cms2} $\dd \va(y_i) \to 0 \quad \text{in} \quad \fs{F}'_c$;
		\item \label{eq:cms3} $\snorm[\fs{F},n]{x_i - y_i} \to 0$ for all $n \in \nn$.
	\end{enumerate}
\end{corollary}
\begin{proof}
	Let $(x_i)_{i \in \nn}$ be a minimizing sequence. This means $\lim_{i \to \infty} \va(x_i) = \inf_{\fs{F}}\va$.
	Define a sequence of positive numbers $(\varepsilon_i)_{i \in \nn}$ as follows
	\[
	\varepsilon_i = \va(x_i) - \inf_{\fs{F}} \va.
	\]
	Since $(x_i)$ is a minimizing sequence, $\lim_{i \to \infty} \varepsilon_i = 0$. Note that if $\va(x_i) = \inf_{\fs{F}}\va$ for some $i$, then $\varepsilon_i = 0$. To ensure $\varepsilon_i > 0$ for the application of Theorem \ref{sec}, we can  modify $\varepsilon_i$ (e.g., set $\varepsilon_i = 1/i$ if $\va(x_i) = \inf_{\fs{F}}\va$, or just use the standard interpretation that $\varepsilon > 0$ in the theorem allows for sequences tending to zero). For the purpose of applying the theorem  we choose a sequence $\varepsilon_i \to 0^+$ such that $\va(x_i) \le \inf_{\fs{F}}\va + \varepsilon_i$. 
	For example, define $\varepsilon_i = \max(\va(x_i) - \inf_{\fs{F}}\va, 1/i)$.

	For each $x_i$ in the minimizing sequence, and for each $\varepsilon_i > 0$ (chosen as described above, ensuring $\varepsilon_i \to 0$), and for each $n \in \nn$ (which determines the index $i$ in conditions (a) and (b) of Theorem \ref{sec}), by virtue of Theorem \ref{sec}, there exists a point $y_i \in \fs{F}$ such that
	\begin{enumerate}[label=\textup{(\alph*)},ref=\textup{\alph*}]
		\item $\snorm[\fs{F},j]{y_i - x_i} \le \dfrac{\va(x_i) - \va(y_i)}{\sqrt{\varepsilon_i}}$ for $j = 1, \ldots, n$;
		\item $\snorm[\fs{F},j]{y_i - x_i} < \sqrt{\varepsilon_i}$ for $j = 1, \ldots, n$;
		\item $\snorm[S]{\dd \va(y_i)} \le K_S \sqrt{\varepsilon_i}$, for all $S \in \mathfrak{S}$ and some $K_S > 0$.
	\end{enumerate}
	From condition (b), since $\varepsilon_i \to 0$, it immediately follows that $\snorm[\fs{F},j]{y_i - x_i} \to 0$ for each $j \in \nn$. This proves condition \eqref{eq:cms3} of the corollary.
	
	From condition (c), since $\varepsilon_i \to 0$, it implies $\snorm[S]{\dd \va(y_i)} \to 0$ for all $S \in \mathfrak{S}$. This means $\dd \va(y_i) \to 0$ in $\fs{F}'_c$, which proves condition \eqref{eq:cms2} of the corollary.
	
 From condition (a), we have $\sqrt{\varepsilon_i} \snorm[\fs{F},j]{y_i - x_i} \le \va(x_i) - \va(y_i)$. Since $\snorm[\fs{F},j]{y_i - x_i} \ge 0$ and $\sqrt{\varepsilon_i} > 0$, this implies $\va(x_i) - \va(y_i) \ge 0$, or $\va(x_i) \ge \va(y_i)$. This proves condition \eqref{eq:cms1}.
	
	Finally, since $\va(x_i) \to \inf_{\fs{F}}\va$ and $\va(y_i) \le \va(x_i)$, and also from $\va(x_0) \le \inf_{\fs{F}}f + \eta$ in Ekeland's principle leads to $\va(z) \le \va(x_0)$, we have $\inf_{\fs{F}}\va \le \va(y_i) \le \va(x_i)$. As $i \to \infty$, $\va(x_i) \to \inf_{\fs{F}}\va$, thus $\va(y_i) \to \inf_{\fs{F}}\va$. Therefore, $(y_i)_{i \in \nn}$ is also a minimizing sequence.
\end{proof}

\begin{corollary}\label{co:minimizing}
	Let  $\va$ be as in Theorem~\ref{sec}. 
	If the $(\mathrm{PS})_m$-condition holds with $m = \inf_{\fs{F}} \va$, then $\va$ attains its minimum at a critical point $ x_0 \in \fs{F} $ with $ \va(x_0) =m $.
\end{corollary}
\begin{proof}
	By Corollary~\ref{cor:almostc}, there exists a minimizing sequence $(y_i)_{i \in \nn}$ such that $\dd \va(y_i) \to 0$ in $\fs{F}'_c$.
	Since $(y_i)$ is a minimizing sequence, we have $\lim_{i \to \infty} \va(y_i) = \inf_{\fs{F}} \va=m$. 
	Thus, for the sequence $(y_i)$, we have $\va(y_i) \to m$ and $\dd \va(y_i) \to 0$. This is  a Palais-Smale sequence at level $m$.
	Since $\va$ satisfies the $(\mathrm{PS})_m$-condition, this sequence $(y_i)$ must have a convergent subsequence. Let $(y_{i_k})$ be such a convergent subsequence, and let $x_0 = \lim_{k \to \infty} y_{i_k}$.
	By the continuity of $\va$ and $\dd \va$ (as $\va$ is a Keller $C_c^1$-functional), we have
	$$ \va(x_0) = \lim_{k \to \infty} \va(y_{i_k}) = m = \inf_{\fs{F}}\va. $$
	And
	$$ \dd \va(x_0) = \lim_{k \to \infty} \dd \va(y_{i_k}) = 0 \quad \text{in} \quad \fs{F}'_c. $$
	Therefore, $x_0$ is a critical point where $\va$ attains its minimum value $m$.
\end{proof}

Now, we prove that the PS-condition for functionals bounded below
implies the coercivity.
\begin{definition}\label{def:coercive}
	Let $ \va : \fs{F} \to \rr $ be a $ C^1 $-functional. We say that $ \va $ is coercive if $ \va(x) \to \infty $ as
	$\snorm[\fs{F},n]{x} \to \infty$ for all $n \in \nn$.
\end{definition}
The following theorem known as Ekeland's $\epsilon$-Principle is a very useful method for finding approximate minima, also called pseudo-minima, which are close to being minima even when no true minimum exists. 

\begin{theorem}[Theorem 1.4.1, \cite{z}]\label{th:ekf}
	Assume that $(\mathsf{M}, \mathbbm{m})$ is a complete metric space.
	Let a functional $\upphi : \mathsf{M} \to (-\infty, \infty]$ be lower semicontinuous, bounded from below, and not identically equal to $\infty$.
	Then, for any $ \epsilon > 0 $ and every point $ m_0 \in \mathsf{M} $ such that $\upphi(m_0) < \infty$, there exists $m_{\epsilon} \in \mathsf{M}$ such that
	\begin{enumerate}
		\item $\upphi(m_{\epsilon}) \le \upphi(m_0) - \epsilon \mathbbm{m} (m_{\epsilon},m_0)$;
		\item $\upphi(m_{\epsilon}) < \upphi(m) + {\epsilon} \mathbbm{m} (m,m_{\epsilon}) \quad \forall m \in \mathsf{M} \setminus \{m_{\epsilon}\}$.
	\end{enumerate}
\end{theorem}
In the following theorem, we apply the translation-invariant metric $\mt$ given by
\[
\mt(x,y) \coloneqq \sum_{n=1}^{\infty} \dfrac{1}{2^n}\dfrac{\snorm[\fs{F},n]{x-y}}{1+\snorm[\fs{F},n]{x-y}},
\]
which generates a topology equivalent to the one defined by the family of seminorms.
\begin{theorem}\label{th:corecive}
	For every \(n \in \nn\), let	$a_n \coloneqq \liminf_{\snorm[\fs{F},n]{x} \to \infty} \varphi(x) < \infty$.
	Then there is a sequence $ (x_i) \subset \fs{F} $ such that as $ \snorm[\fs{F},n]{x_i} \to \infty $, we have
	$$\varphi (x_i) \to a_n,\,  \va'(x_i) \to 0 .$$
\end{theorem} 

\begin{proof}
	For every \(n \in \nn\), define the function 
	\begin{equation}\label{corecive:mr}
		\bl{m}_n(r) \coloneqq \inf_{\snorm[\fs{F},n]{x} \geq r} \varphi(x) \quad r \in \rr_{+}.
	\end{equation}
	The function $ \bl{m}_n(r) $ is a non-decreasing, and 
	\begin{equation}\label{corecive:ee}
		\lim_{r \rightarrow \infty} \bl{m}_n(r) = a_n.
	\end{equation} 
	By~\eqref{corecive:ee} for each $  \epsilon >0 $
	there exists $ r_1   $ such that for all $ r\geq r_1 $, we get
	\begin{equation} \label{corecive:27}
		a_n -\epsilon^2 \leq \bl{m}_n(r).
	\end{equation}
	Now, for a fixed $ \epsilon >0 $ choose a number $ r_2 $ such that
	\begin{equation}\label{corecive:46}
		r_2 \geq \max \{ r_1,2\epsilon\}.
	\end{equation}
	By our assumption, we can fix some $ x_0 $ with $  \snorm[\fs{F},n]{x_0} \geq 2r_2$ such that
	\begin{equation}\label{corecive:42}
		\varphi(x_0) <  a_n + \epsilon^2.
	\end{equation} 
	Let $ \mathbf{F} = \Set{x \in \fs{F}}{\snorm[\fs{F},n]{x} \geq r_2} $. This is a closed subset of $ \fs{F} $, so it is a complete metric space with the metric induced from \(\fs{F}\). Moreover, $ \varphi $ is lower semicontinuous on $ \fs{F} $ and so on $ \mathbf{F} $. Also, by
	\eqref{corecive:mr},\eqref{corecive:27} and \eqref{corecive:46} we obtain
	\[ \varphi(u)  \geq  \bl{m}(\snorm[\fs{F},n]{u}) \geq a_n - \epsilon^2 \quad \forall u \in \fs{F} \, \text{ with } \, \snorm[\fs{F},n]{u} \geq r_2.
	\]
	So, $ \varphi $ is lower bounded, and therefore, all assumptions of Theorem~\ref{th:ekf} are fulfilled for  
	$\mathbf{F} $. Thus, there is $ x_{\varepsilon} \in \mathbf{F} $ such that 
	\begin{equation} \label{corecive:lklk}
		\varphi(x_{\epsilon}) <  \varphi(x) + \epsilon \mt(x_{\epsilon},x), \quad \forall x \in \mathbf{F}\setminus \set{x_{\epsilon}},
	\end{equation}
	and, as well as
	\begin{equation}\label{corecive:24}
		\varphi (x_{\epsilon}) \leq \varphi(x_0) - \epsilon  \mt(x_{\epsilon},x).
	\end{equation}
	It follows from \eqref{corecive:mr},\eqref{corecive:27},\eqref{corecive:46},\eqref{corecive:24} and \eqref{corecive:42} that
	\begin{equation*}
		a_n - \epsilon^2 \leq \bl{m}_n(r_2) \leq \varphi(x_{\epsilon}) \leq \varphi(z_0)-\epsilon 
		\mt(x_{\epsilon},x)   \leq a_n + \epsilon^2 - \mt(x_{\epsilon},x).
	\end{equation*}
	Hence, we obtain $\mt(x_{\epsilon},x) \leq 2 \epsilon$. 	
	Thereby,  \eqref{corecive:46} implies
	\begin{equation*}
		\mt(x_{\epsilon},\zero{\fs{F}})   \geq \mt(x_{0},\zero{\fs{F}})- \mt(x_{\epsilon},x) \geq 2r_2 - 2 \epsilon \geq r_2.
	\end{equation*}
Here, \(\zero{\fs{F}}\) is the origin of \(\fs{F}\).
	Whence  
	$ x_{\epsilon} $ is an interior point of $ \mathbf{F} $, and by Theorem~\ref{sec} we obtain
	\begin{equation}\label{eq:mi}
		\snorm[S]{\va'(x_{\epsilon})}< \epsilon \quad  \forall S \in \bb.
	\end{equation}
	For each $i \in \nn$, taking $\varepsilon = 1/i$, we get a sequence \(x_i\) such that by \eqref{corecive:42}
	and \eqref{corecive:24}, it satisfies \(\varphi (x_i) \to a_n\). Moreover, by  \eqref{eq:mi}, it satisfies \(\va'(x_i) \to 0\).
\end{proof}
\begin{corollary}\label{cor:coercive}
	If the functional $ \va \colon \fs{F} \to \rr $ is bounded below and satisfies the $\mathrm{PS}$-condition for every $ m \in \rr $, then $ \va $ is coercive.
\end{corollary}

\begin{proof}

	If $\va$ is not coercive, then the condition $\va(x) \to \infty$ as $\snorm[\fs{F},n]{x} \to \infty$ for all $n \in \nn$ fails. This implies that there is at least one $n_0 \in \nn$ for which $\liminf_{\snorm[\fs{F},n_0]{x} \to \infty} \va(x) = a_{n_0} < \infty$.
	
	In this case, Theorem~\ref{th:corecive} implies  that there exists a sequence $(x_i) \subset \fs{F}$ such that for this specific $n_0$, $\snorm[\fs{F},n_0]{x_i} \to \infty$, and simultaneously $\va(x_i) \to a_{n_0}$ and $\dd \va(x_i) \to 0$ in $\fs{F}'_c$.

	Since $\va$ satisfies the $\mathrm{PS}$-condition at level $a_{n_0}$, it follows that the sequence $(x_i)$ must possess a convergent subsequence.
	However, a convergent sequence must be bounded. If a subsequence of $(x_i)$ converges, it must be bounded with respect to all seminorms. This contradicts that $\snorm[\fs{F},n_0]{x_i} \to \infty$ for this specific $n_0$.
\end{proof}


\begin{thebibliography}{99}

\bibitem{eftekharinasab}
K.~Eftekharinasab and I.~Lastivka, \emph{A Lusternik-Schnirelmann type theorem for \(C^1\)-Fréchet manifolds}, Journal of the Indian Mathematical Society \textbf{88} (2021), no.~3--4, 309--322.

\bibitem{z}
C.~Zălinescu, \emph{Convex Analysis in General Vector Spaces}, World Scientific, 2002.

\bibitem{ke}
H.~Keller, \emph{Differential Calculus in Locally Convex Spaces}, Springer-Verlag, Berlin, 1974.

\bibitem{gl}
H.~Glöckner, \emph{Aspects of differential calculus related to infinite-dimensional vector bundles and Poisson vector spaces}, Axioms \textbf{11} (2022), no.~5, 221.

\bibitem{mm}
J.~Mawhin and M.~Willem, \emph{Origin and evolution of the Palais-Smale condition in critical point theory}, J.~Fixed Point Theory Appl. \textbf{2} (2010), 265--290.

\bibitem{ps}
R.~S.~Palais and S.~Smale, \emph{A generalized Morse theory}, Bull.~Amer.~Math.~Soc. \textbf{70} (1964), 165--172.

\bibitem{ekel}
J.~H.~Qiu, \emph{Ekeland's variational principle in Fréchet spaces and the density of extremal points}, Studia Math. \textbf{1} (2005), 81--94.

\bibitem{tw}
T.~Wurzbacher, \emph{Fermionic second quantization and the geometry of the restricted Grassmannian}, in \emph{Infinite Dimensional Kähler Manifolds (Oberwolfach, 1995)}, DMV Sem., vol.~31, Birkhäuser, Basel, 2001, pp.~287--375.

		
	\end{thebibliography}
\end{document}